\documentclass{article}
\usepackage{amsmath,amsfonts,amsthm,amssymb}
\usepackage{graphicx}

\newcommand{\Z}{\mathbb{Z}}

\newcommand{\cbrt}{ \mathop \mathrm{cbrt} \nolimits}
\newcommand{\ph}{ \mathop \mathrm{ph} \nolimits}
\newcommand{\Ln}{ \mathop \mathrm{Ln} \nolimits}

\newlength{\funcwd}\newlength{\funcht}

\newcommand{\invItal}[1]{
\settowidth{\funcwd}{#1}\settoheight{\funcht}{#1}%
\raisebox{.9\funcht}{%
 \makebox[0pt][l]{\hspace{.6\funcwd}$\scriptstyle\vee$}}#1} 

\begin{document}

\title{Branch Cuts and Riemann Surfaces}
\author{ David Jeffrey \\
Department of Mathematics, The University of Western Ontario \\
London, Ontario, Canada
}

\date{}

\maketitle

\begin{abstract}
  The plotting of Riemann surfaces by computational software is discussed. The link between the branches of a multi-valued
  function $g(z)$, defined on the range of $g(z)$, and a Riemann surface, defined on the domain of $g(z)$, is emphasized.
  The connection between the two is clarified by defining the \textit{charisma} of the argument $z$ to the function.
  This approach is used to plot several surfaces using Maple.
\end{abstract}
\section{Introduction}
Riemann surfaces of complex functions have been visualized by a number of authors and computer systems, and
often make striking visual images. See \cite{Trott} for an impressive collection. The traditional way that textbooks
introduced Riemann surfaces is through an informal discussion of layered sheets, which are cut and glued together to form a continuous (usually self-penetrating) surface \cite{Wegert}.
Such an approach does not lend itself to using software to create the surfaces.
This was discussed in \cite{Corless1998Riem}, prompted by a demonstration program in Matlab.
The methods used there are extended here to strengthen the link between the range and the domain.

Riemann surfaces are used to understand multi-valued functions, which are often the inverses of `proper' single-valued functions. Examples include the inverses of the exponential function (namely logarithm), the trigonometric functions (arc-trigs),
and the integral powers ($n$th-roots).
In each case, there is a function $z=f(w)$ which is single valued, and an inverse $w=\invItal{f}(z)$ which is multivalued.
For the multivalued function, the question is how to separate and access the various elements in the set of multiple values.
Two possibilities present themselves: the separation is made either in the range of the function or in its domain.

\section{Labelling the range}\label{sec:range}
We begin with a common example of a multivalued function: the logarithm.
If $ z=e^w$, then $w =\ln z $.
Since, for $n\in \Z$, $ z=e^{w +2n\pi i} = e^w$, for any complex number $z$, there are an infinite number of values for the logarithm.
The standard treatment defines two logarithm functions.
Following the notation of the DLMF \cite{AS, DLMF}, we write $\Ln z$ to
represent the \textit{general logarithm function}, which stands for the infinite collection of values, and $\ln z$ for the
\textit{principal value} or \textit{principal branch}. The principal branch being defined by $-\pi< \Im(\ln z) \le \pi$.
The general function notation, however, is frustratingly vague, and leads to statements such as
\begin{equation}\label{eq:genlog}
  \Ln z = \ln z + 2k\pi i\ ,
\end{equation}
where the equation has one variable $z$ on the left and two variables $z,n$ on the right.

\begin{figure}
  \centering
  \includegraphics[width=5cm]{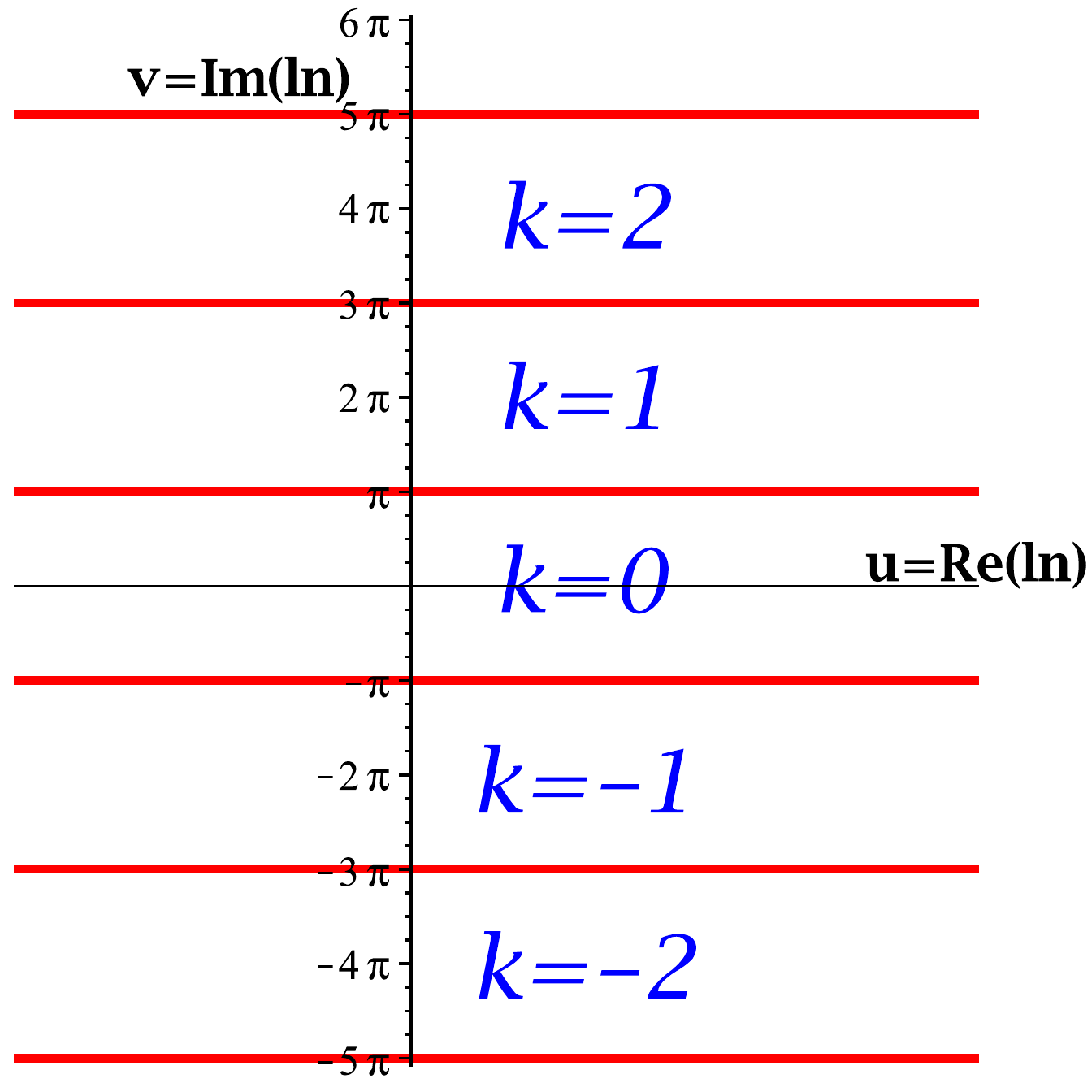}
  \caption{The range of the logarithm partitioned into branches. The branches are labelled with the integer $k$, and the index is added
  to the function notation: $\ln_k z$.
  Thus, to stipulate a logarithm value with imaginary part between $\pi$ and $3\pi$, one writes
  $\ln_1 z$.}\label{fig:lnbranch}
\end{figure}

In order to obtain a more precise definition of the logarithm value, the notation
\begin{equation}\label{eq:definelnk}
  \ln_k z=\ln z +2k\pi i \ ,
\end{equation}
was introduced in \cite{unwindW}.
Geometrically, this coincides with the range of logarithm being partitioned into branches, and each branch labelled.
See figure \ref{fig:lnbranch}. We note $\ln_0 z$ also denotes the principal branch.
For each separate branch $\ln_k z$, the domain consists of the whole complex plane, with a line of discontinuity, called a branch cut, along the negative real axis. In thinking of branches, it is sometimes helpful to disregard the connections between the branches and think of each branch as a function separate in its own right. For example, the equation $\ln x=2$ has the solution $x=e^2$, but the equation $\ln_1 x=2$ has no solution because $2$ does not lie in the range of $\ln_1 x$.

\begin{figure}
  \centering
  \includegraphics[scale=0.35]{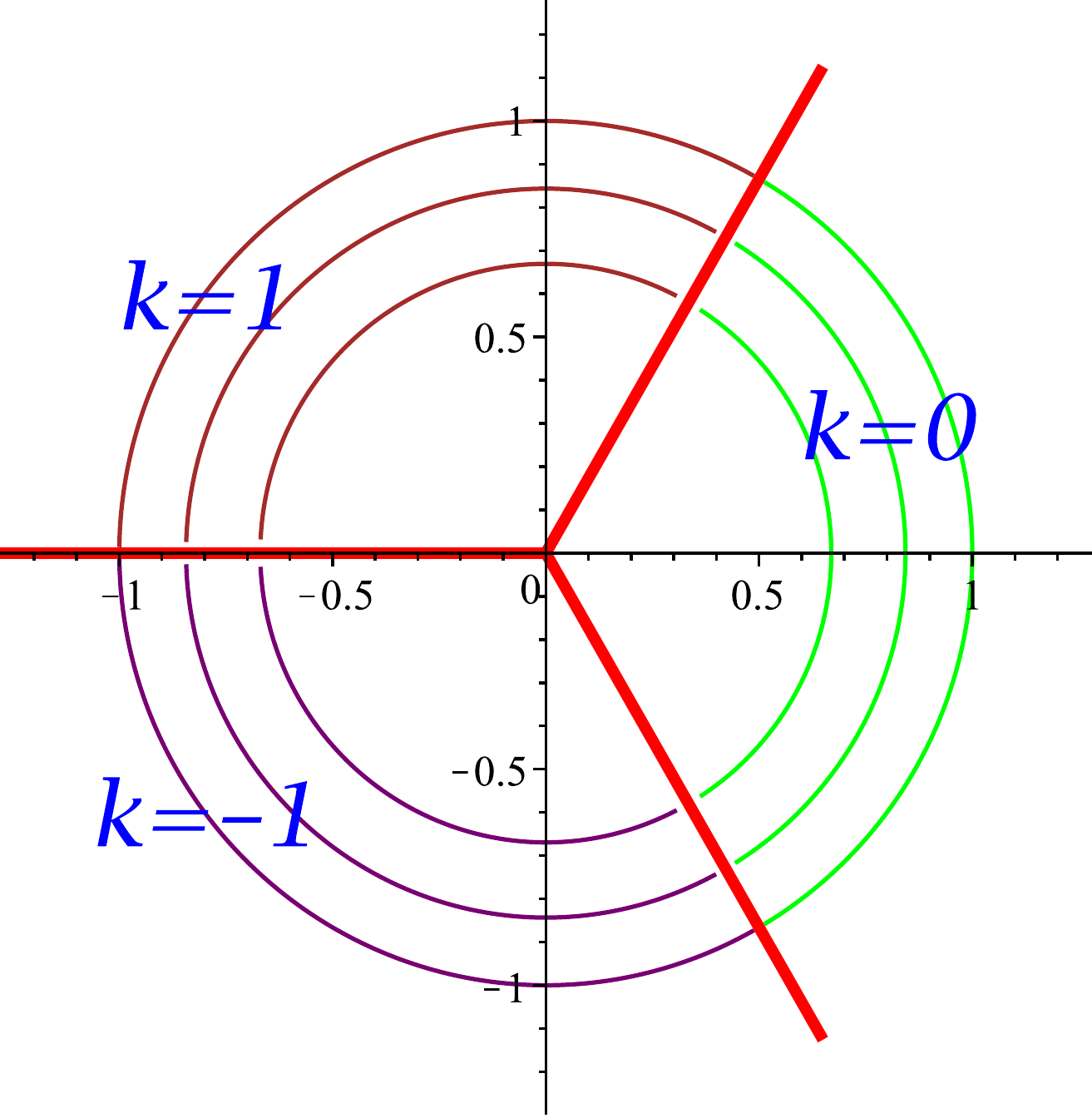}
  \caption{The range of the cube root partitioned into branches. As with logarithm, the branches are labelled with an integer $k$. In this case there are only 3 branches, separated by the boundaries $r e^{\pm i\pi/3}$ and the negative real axis.}\label{fig:cbrtbranch}
\end{figure}

A second example is provided by the cube-root function. We have $z=w^3$, implying the inverse $w=z^{1/3}$.
In the complex plane, the cube root always has 3 values. To denote the 3 cube roots, we need a notation which gives us space for a label.
Since the notations $z^{1/3}$ and $\sqrt[3]z$ look clumsy with additional labels, we create a cube-root function name:
$\cbrt_k z$.
The range of $\cbrt$ is partitioned as shown in figure \ref{fig:cbrtbranch}, where the branches are defined using the complex
phase\footnote{We use the name complex phase, rather than complex argument, because the word argument will used for function argument.}.
With $\ph z$ denoting the principal complex phase $-\pi<\ph z \le\pi$, we set the principal branch by the requirement
 $-\pi/3<\ph(z^{1/3})\le \pi/3$. If we denote the primitive root of unity by $\omega=e^{2\pi i/3}$, and $\theta=\ph z$,
 then the 3 branches are defined by
\begin{equation}\label{eq:cbrtbranchlabel}
 \cbrt_k(z) = \cbrt_k(re^{i\theta})
              =\begin{cases}
                     \phantom{\omega}r^{1/3} e^{i\theta/3}\ , & k=0\ ,\\
                     \omega r^{1/3} e^{i\theta/3}\ , & k=1\ , \\
                     \overline \omega r^{1/3}e^{i\theta/3} \ , & k=-1\ .
               \end{cases}
\end{equation}
\noindent
For each separate branch of $\cbrt_k(z)$, the
domain is the whole of the complex plane, with a line of discontinuity (the branch cut) along the negative real axis.

We note as an instructive special case the values of $(-8)^{1/3}$. The three cube roots are $-2,1+\sqrt3\,i,1-\sqrt3\, i$.
Some readers might be surprised to learn that the principal-branch value of $(-8)^{1/3}$ is not $-2$.
The principal branch value is $\cbrt_0(-8) = 1+\sqrt3\, i$ and the real-value root is branch 1: $\cbrt_1(-8)=-2$.
As seen in figure \ref{fig:cbrtbranch}, the negative real values of the function do not lie in the principal branch.

We can note here that the principal-branch values of both logarithm and cube root are the unique values returned by all major scientific
software, such as Matlab, Maple and Mathematica. This also applies to other multi-valued functions, such as $n$th roots and inverse trigonometric functions.

\section{Labelling the domain}
The main object of this paper is to show how the labelling ideas of the previous section can be transferred from the range of a complex function
to its domain, and by doing so we obtain a new perspective on plots of Riemann surfaces.
A Riemann surface for a multivalued function is built on its domain. A typical description of the construction process
talks about cutting and joining sheets of the domain.
For example, here is a description of a Riemann surface for the square-root function \cite{brown}.
\begin{quote}
  A Riemann surface for $z^{1/2}$ is obtained by replacing the $z$ plane with a surface made up of two sheets $R_0$ and $R_1$, each cut along the positive real axis with $R_1$ placed in front of $R_0$. The lower edge of the slit in $R_0$ is joined to the upper edge of the slit in $R_1$, and the lower edge of the slit in $R_1$ is joined to the upper edge of the slit in $R_0$.
\end{quote}


A Riemann surface is erected over the domain of a multivalued function, that is, given $w=\invItal{f}(z)$, the surface represents values of $z$ rather than values of $w$.
The value of $z$ by itself is not sufficient to determine the value of $\invItal{f}(z)$,
and therefore there must be a property possessed by a particular value of $z$ that decides the value of the function, in conjunction with the complex value of $z$. This property is not at present visible.
We shall call this property \textit{charisma}. A variable $z$ with charisma will be denoted $z_{\cal C}$\,, while we decide what it is.

We aim to define charisma as a numerical value, which can then serve as an ordinate on an axis perpendicular to the complex plane.
We set up 3 axes: real and imaginary axes for locating the complex value of $z$ together with an orthogonal axis representing charisma.
The first example will be the cube-root function described above.

\subsection{Charisma for cube root}
We take as a first example the indexed cube root given in \eqref{eq:cbrtbranchlabel}.
We explore four possible choices for the charisma of this function.

\subsubsection{Charisma as branch index}
\begin{figure}
  \centering
  \includegraphics[scale=0.5]{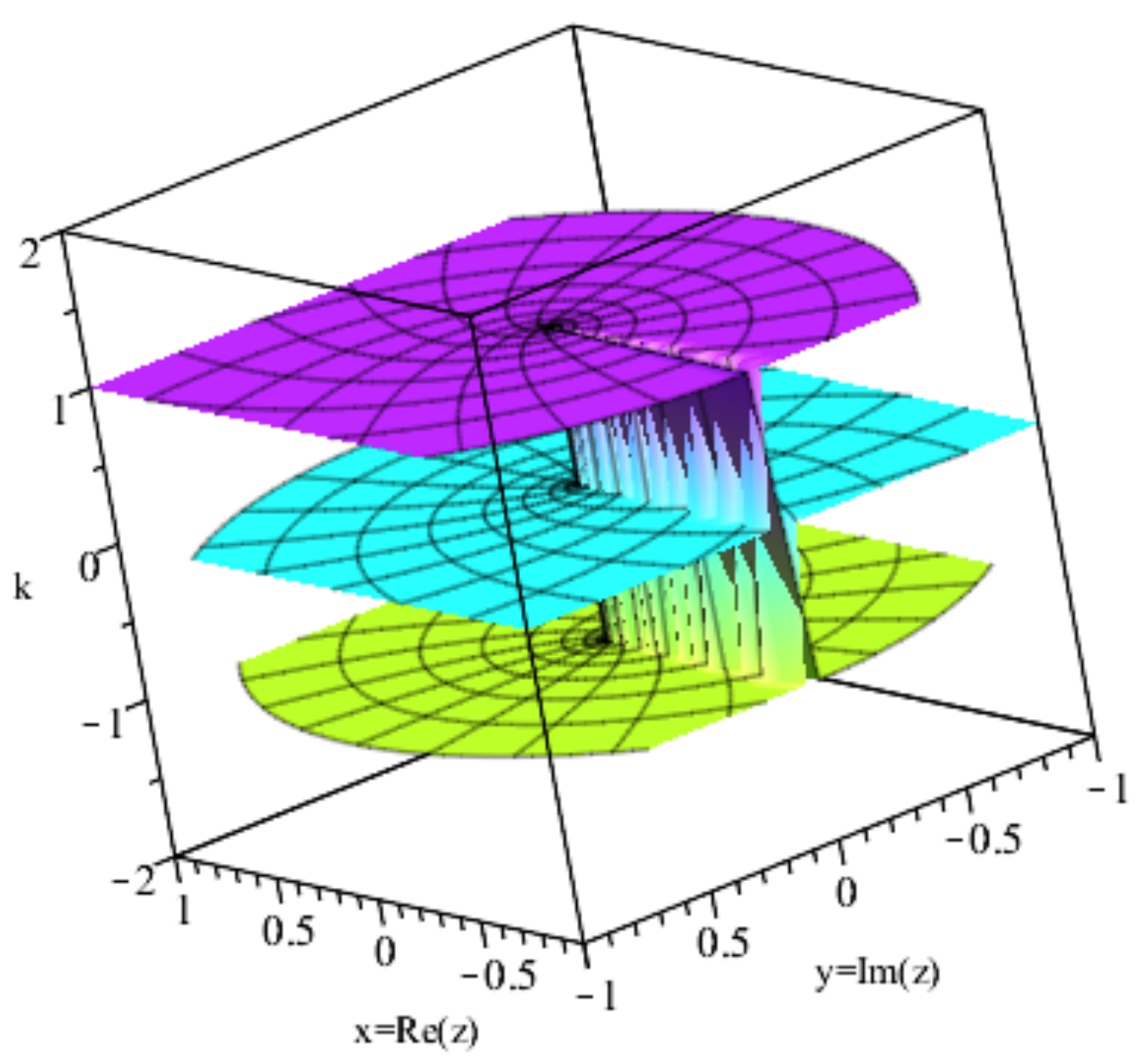} \quad
\raisebox{0cm}{\includegraphics[scale=0.3]{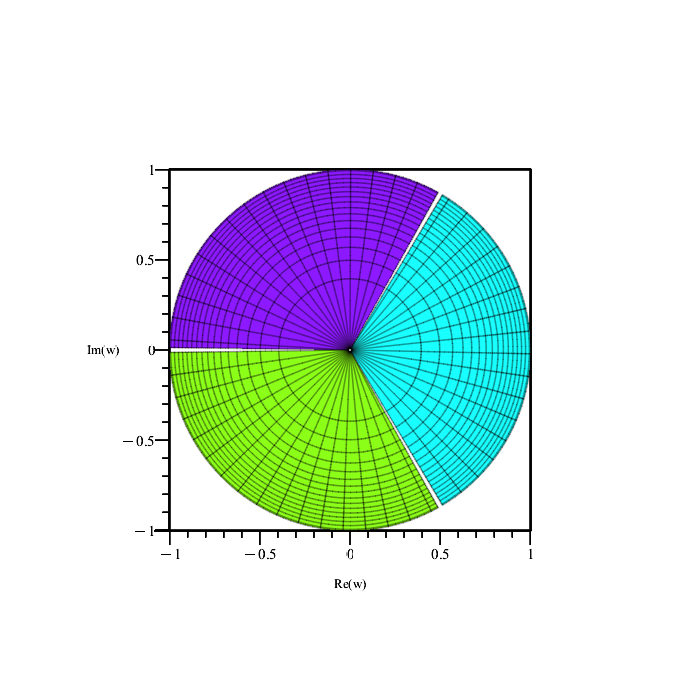}}
  \caption{On the left we see the domain of $z^{1/3}$ as a Riemann surface; on the right, the range of $z^{1/3}$. The sheets of the Riemann surface are separated by labelling them with the branch index. The branches of $\cbrt_k(z)$ are coloured in the same way as the sheets of the Riemann surface. The top sheet corresponds to $k=1$ meaning all the values of $z$ on this sheet map to $\cbrt_1(z)$ as shown by the corresponding colour in the range. Similarly for the other two branches. One can further note that the sheets are connected discontinuously, reflecting the discontinuous change in branch number as one moves around in the range.}
  \label{fig:charlabel}
\end{figure}

We have seen the branch label $k$ used to define values in the range of $\cbrt_k(z)$, so we begin by trying the assignment
\[{\cal{C}} =k\ .\]
With the charisma chosen to be the index $k$, we plot the corresponding surfaces as follows. We generate an array of points in the complex plane.
Then we create a 3-element plot structure consisting of the real and imaginary parts of a complex number, together with the charisma.
This structure is applied to the array of complex values three times, once for each of the three values of the charisma.
The complete array is handed to a three-dimension plotting command in Maple.
The pictures we get are shown in figure \ref{fig:charlabel}.

The effect of the charisma is simply to lift the flat complex planes so that the three planes are stacked and spaced.
The vertical walls in the plot are Maple's way of showing that the planes are connected discontinuously.
We can understand this plot by imagining a point\footnote{Or more picturesquely an ant. }  in the \textit{range} of the function.
The point now circles the origin in the range. Each time the point crosses a branch boundary, the charisma jumps discontinuously with the
branch number. This discontinuous change is reflected in the domain by a jump from one sheet to another.

\subsubsection{A continuous charisma} \label{sec:CharismaPhase}
Our first choice has an unsatisfactory feature. If we follow our point in the range of $z^{1/3}$,
the complex values it samples change continuously. The discontinuity arises solely from the branch label being discontinuous.
It is the same as driving across the boundary between two provinces or states. The road is continuous, the land is continuous, but suddenly everyone is speaking French and selling fireworks\footnote{At least that is what I see when driving from Toronto to Montreal.}.

\begin{figure}
\includegraphics[scale=0.3]{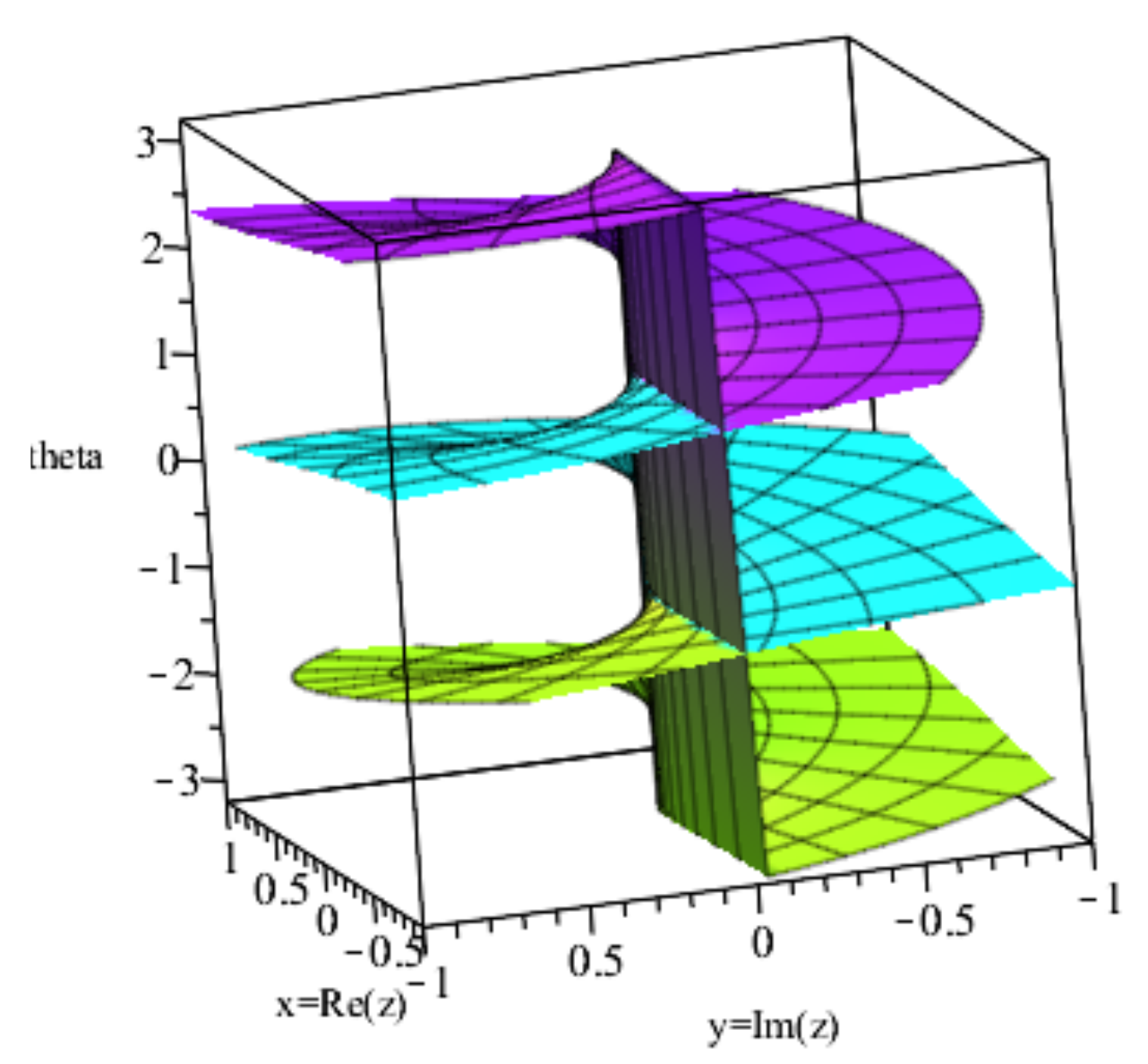} \quad
\raisebox{0cm}{\includegraphics[scale=0.2]{cbrt.png}}
\caption{The Riemann surface obtained by using the phase of the cube root as charisma. The phase has a minimum value of $-\pi$ and increases monotonically to $-\pi/3$ while the branch $k=-1$ covers its domain, which is the whole of the complex plane. This is the green surface in the domain. The behaviour is repeated for the other two branches. The surfaces meet smoothly where they join. Maple has joined the end of the surface at $\mathcal{C}=\pi$ discontinuously to $\mathcal{C}=-\pi$.}\label{fig:charismaphase}
\end{figure}

In the previous section, we see that the simplest choice of charisma misses an important property of the
cube root, namely that between the branches there is a continuous transition. We wish to find a charisma that
captures this behaviour.
When searching for a measure of charisma, it is important\footnote{Too strong a word? Well, at least helpful.} not to fixate on the domain of $z$, even though the domain is where the surface will be located. We are trying to represent the multi-valuedness of the function, and that is defined in the range.
It was remarked that a point circling the origin in the range would see continuous behaviour.
Following the implications of this observation suggests that the phase of the cube root could be a better quantity to use for charisma.
Thus, given $w=\cbrt_k(z)$, we set $\mathcal{C}=\ph(w)$. The phase is computed from the value of $w$, the range of $z^{1/3}$, not from the value of $z$. This choice results in the Riemann surface shown in figure \ref{fig:charismaphase}.
This removes the jumps seen in the previous section, but is still unsatisfactory in that
in the range of $z^{1/3}$
the branch $k=1$ joins smoothly to branch $k=-1$,
but the Riemann surface in figure \ref{fig:cbrtbranch} joins discontinuously.
Therefore, we try now a third choice which will give us smooth joins of all of the branch transitions.
\subsubsection{A continuous and periodic charisma}
\begin{figure}
\begin{center}
\includegraphics[scale=0.5]{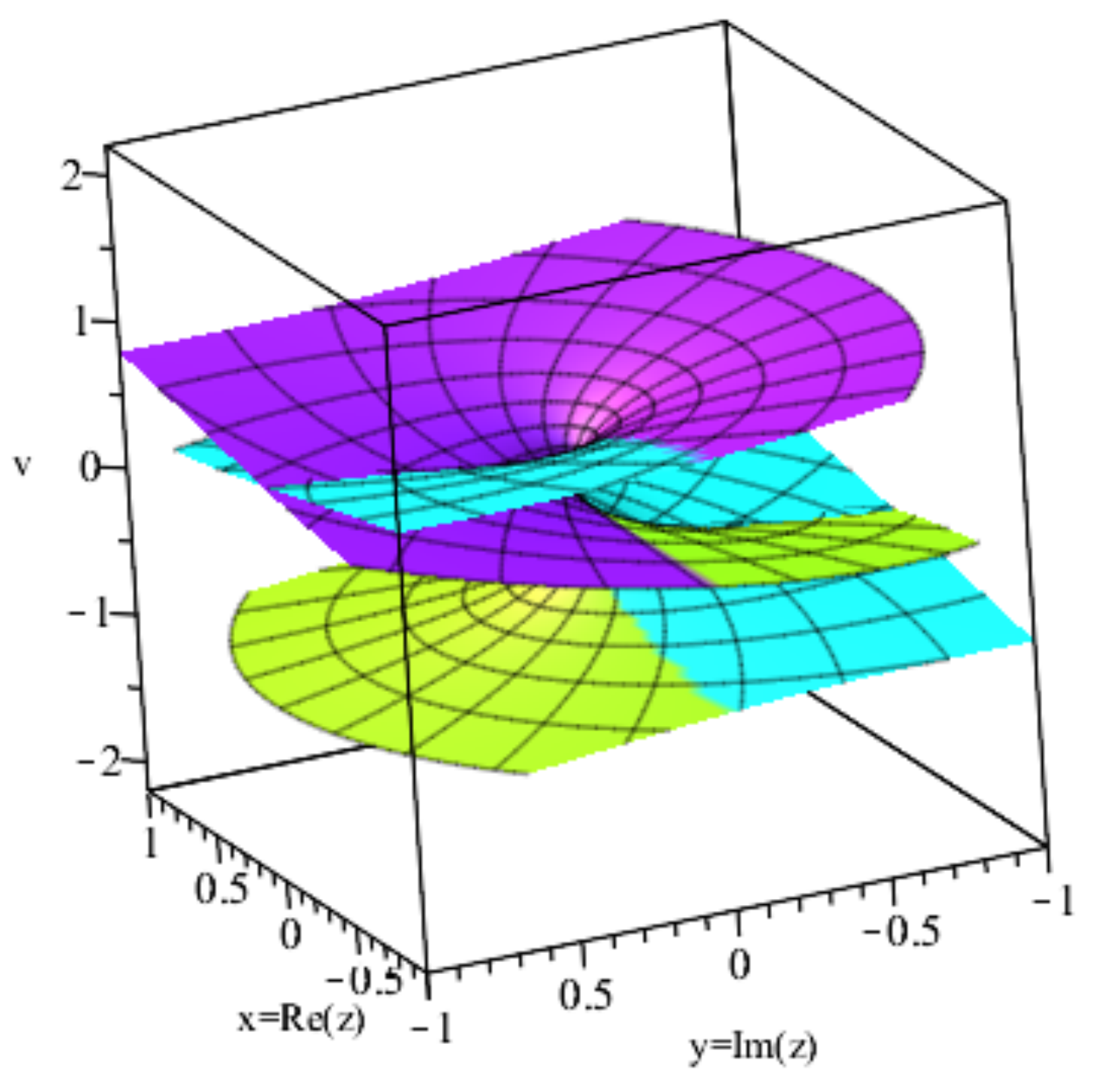}
\end{center}
\caption{A periodic Riemann surface for the cube-root function.
The charisma is based on $\sin \theta$ where $\theta=\ph(z^{1/3})$.
This can equivalently be based on the imaginary part of $z^{1/3}$.
Since the Riemann surface is plotted over the domain $z$, while being coloured based on the range $z^{1/3}$,
the colours change along the negative real axis in the domain where $\mathcal{C}=-\sqrt3/2,0,\sqrt3/2$. }
\label{fig:cbrtsin}
\end{figure}

In the range of $z^{1/3}$ we see that the behaviour is essentially periodic, in that for $x<0$
\begin{equation}\label{eq:periodiccbrt}
\lim_{y\downarrow 0} \cbrt_1(x+iy)=\lim_{y\uparrow 0} \cbrt_{-1}(x+iy) \ .
\end{equation}
A function that has these properties is $\mathcal{C}=\sin\theta\propto \Im(w)$.
In figure \ref{fig:cbrtsin}, the Riemann surface produced with this choice is shown.
The surface is continuous and smooth everywhere.
It now intersects itself in two places: branch $k=-1$ intersects $k=0$ when $\mathcal{C}=-1/2$
along the negative imaginary axis, and $k=0$ intersects $k=1$ at $\mathcal{C}=1/2$ along the positive imaginary axis.
The surfaces change colour (branch) along the negative real axis at
$\mathcal{C}=-\frac{\sqrt3}{2},0,\frac{\sqrt3}{2}$.
\subsubsection{An alternative continuous and periodic surface}
\begin{figure}
\begin{center}
\includegraphics[scale=0.5]{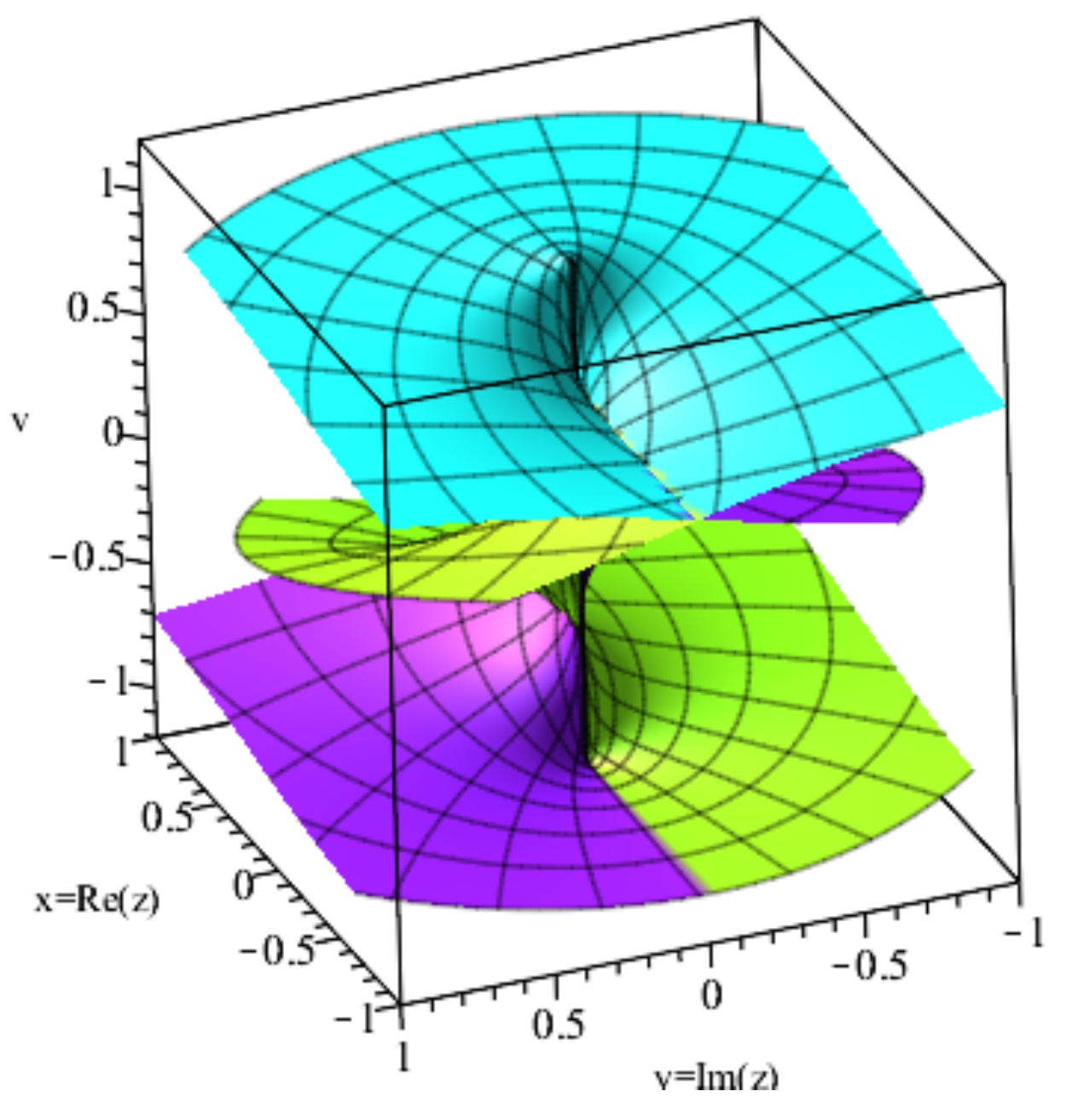}
\end{center}
\caption{The Riemann surface for cube root, using $\cos\theta$ for the charisma. The surface is the same form as the one using the sine function, but the partition between the branches has changed. }\label{fig:cbrtcos}
\end{figure}

It is also possible to use the charisma $\mathcal{C}=\cos\theta$. The result is shown in figure \ref{fig:cbrtcos}. Note that the principal branch now is prominent at the top of the surface.
%
%
%
\subsection{And lastly, logarithm}

\begin{figure}
  \centering
  \includegraphics[scale=0.35]{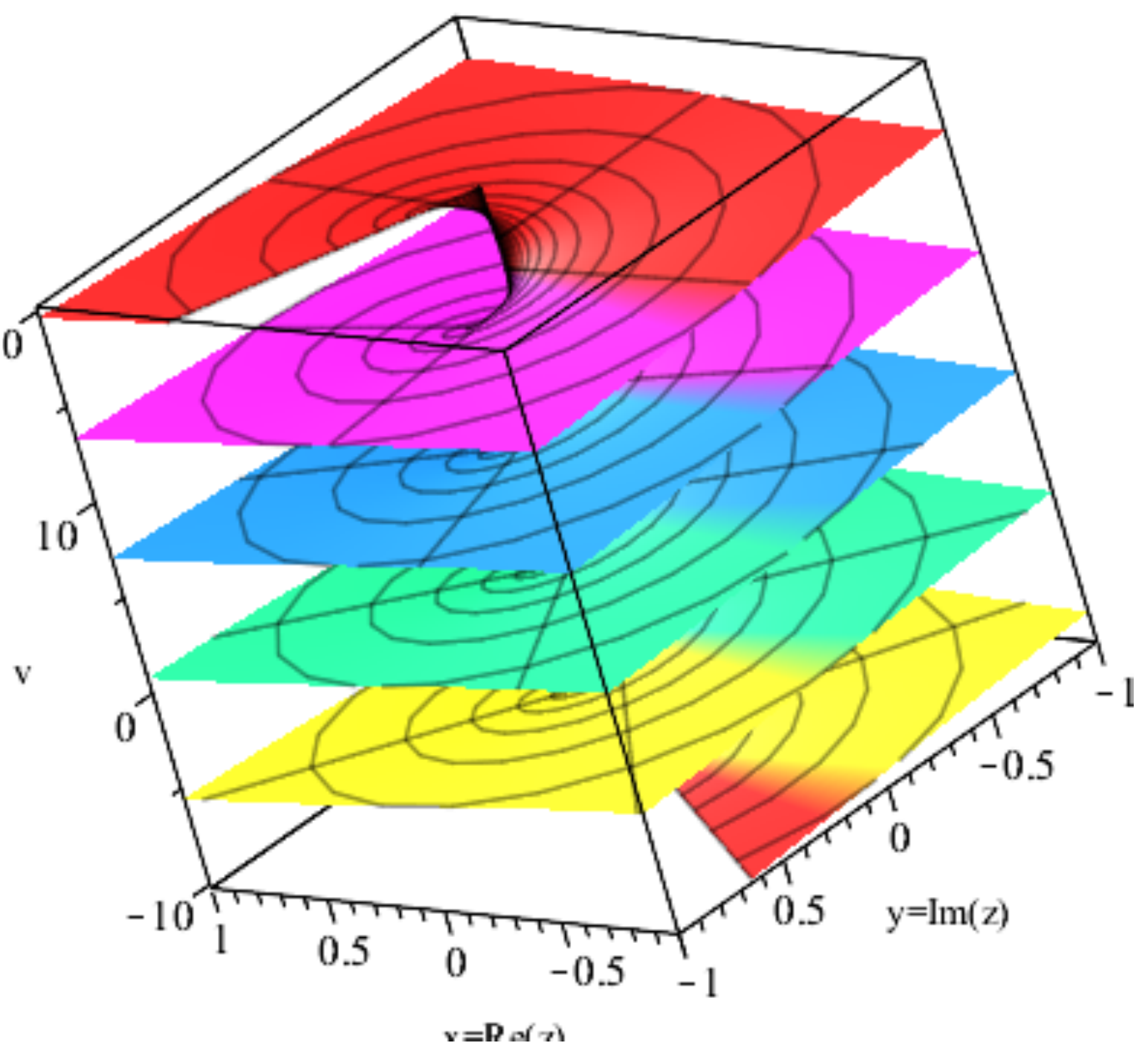}
  \caption{A section of the infinite Riemann surface for logarithm, coloured according to the branch index.}
  \label{fig:logRiemann}
\end{figure}

In \cite{Corless1998Riem}, a contrast was drawn between attempts to plot a Riemann surface for the logarithm.
One attempt was shown to be unsatisfactory. By emphasizing the need to think in terms of the range of a function,
we can see immediately why this is. 
The charisma was chosen to be the real part of the logarithm, 
but a glance at the range of log shows immediately that motion parallel to the real axis 
never leaves the branch of log in which the point started.
In order to obtain a surface for logarithm, we must move parallel to the imaginary axis 
in order to cross from one branch to the next. This was shown to be the correct choice.
We could try using the branch index of $\ln_k z$ as we did for cube root, 
but again that would result in flat segments separated by jumps. 
As shown in \cite{Corless1988Comp}, a simple charisma for logarithm is its imaginary part. 
The result is shown in figure \ref{fig:logRiemann}, where the colouring again changes with the branch.

\section{Conclusions}
Many treatments of multivalued functions concentrate on the domain of the function, and focus on branch cuts in the domain.
See for example, \cite{AS,DLMF,K}. This article hopes to demonstrate thatthinking in the range as well as the domain makes understanding easier.


\end{document}